\documentclass[12pt]{amsart}
\usepackage{amscd,amsmath,amssymb,amsfonts}
\usepackage[all]{xy}
\usepackage{hyperref}
\usepackage{url}
\usepackage{stmaryrd}
\usepackage{color}

\theoremstyle{plain}
\newtheorem{thm}{Theorem}
\newtheorem{lem}[thm]{Lemma}
\newtheorem{cor}[thm]{Corollary}
\newtheorem{prop}[thm]{Proposition}

\newtheorem{conj}[thm]{Conjecture}
\newtheorem{claim}[thm]{Claim}
\theoremstyle{definition}
\newtheorem{defn}[thm]{Definition}
\newtheorem{ex}[thm]{Example}

\newtheorem{rmk}[thm]{Remark}

\numberwithin{thm}{section}

\newcommand{\ga}[2]{\begin{gather}\label{#1}#2 \end{gather}}

\newcommand{\Hom}{{\rm Hom}}

\newcommand{\Spec}{{\rm Spec \,}}

\newcommand{\sG}{{\mathcal G}}

\newcommand{\sM}{{\mathcal M}}
\newcommand{\sN}{{\mathcal N}}
\newcommand{\sO}{{\mathcal O}}

\newcommand{\sQ}{{\mathcal Q}}

\newcommand{\sS}{{\mathcal S}}
\newcommand{\sT}{{\mathcal T}}

\newcommand{\sW}{{\mathcal W}}
\newcommand{\sX}{{\mathcal X}}

\newcommand{\sZ}{{\mathcal Z}}
\newcommand{\A}{{\mathbb A}}

\newcommand{\C}{{\mathbb C}}

\newcommand{\F}{{\mathbb F}}
\newcommand{\G}{{\mathbb G}}

\renewcommand{\P}{{\mathbb P}}

\newcommand{\Z}{{\mathbb Z}}

\newcommand{\et}{{\rm \acute{e}t}}

\newcommand{\Ch}{{\rm Ch}}

\DeclareMathOperator{\GL}{GL}

\begin{document}
\title[Quasi-unipotent monodromy]{ Local systems with quasi-unipotent monodromy at infinity are dense  }
\author{H\'el\`ene Esnault \and Moritz Kerz}
\address{Freie Universit\"at Berlin, Arnimallee 3, 14195, Berlin,  Germany}
\email{esnault@math.fu-berlin.de}
\address{   Fakult\"at f\"ur Mathematik \\
Universit\"at Regensburg \\
93040 Regensburg, Germany}
\email{moritz.kerz@mathematik.uni-regensburg.de}
\thanks{The first  author thanks the Institute for Advanced Study  where this work was initiated.  The second author is supported the SFB 1085 Higher Invariants, Universit\"at Regensburg. Both authors thank the Mathematisches Forschungsinstitut in Oberwolfach where this work was done while they enjoyed the program 'Research in Pair'.}

\begin{abstract} We show that complex local systems with quasi-unipotent monodromy at infinity over a normal complex variety are Zariski dense in their moduli.  \end{abstract}

\maketitle

\section{Introduction}\label{sec:intro}

Let $G$ be a linear algebraic group over the complex numbers.
In this short note we study $G$-representations of the topological fundamental group
$\pi:= \pi_1(X(\C ),x) $ of a  normal complex variety $X$ which are quasi-unipotent with
respect to the monodromy at infinity.
As $\pi$ is finitely generated, the set of group homomorphisms (called
\emph{$G$-representations}) $\rho\colon \pi \to G(\C )$ is in a canonical way
the set of complex points of an affine complex variety
$\Ch_{G, \C}^\Box(\pi) $, the so called \emph{framed character variety}, see Section~\ref{sec:Grep}.

Our main result is motivated by the following conjecture about the density of representations of
geometric origin. We fix an embedding of linear algebraic groups
\[
G \stackrel{\iota}{\hookrightarrow} \GL_{r,\C }.
\]
We say that a $G$-representation $\rho\colon \pi \to G(\C)$ is of  geometric origin   if there is a smooth projective morphism $f\colon Y\to U$, where
$j\colon U\hookrightarrow X$ is a dense open subvariety such that the semi-simplification of the representation \[\iota \circ \rho \circ   j_* \colon \pi_1(U(\C)) \to
\GL_{r,\C}\] gives rise to a linear  local system  which is a direct summand of the local system
$ \ \oplus_i R^if_* \C$ on $U$.    As  linear systems of geometric origin are compatible with tensor products, direct sums and duals, 
 \cite[Thm.4.14]{Mil17} implies that the notion of  $\rho$  being of geometric origin does not depend on the choice of the embedding $\iota$.

\begin{conj}[Density]\label{dens:conj}
 The set of $G$-representations of geometric origin  is Zariski dense in $\Ch_{G, \C}^\Box(\pi) $. {\footnote{Aaron Landesman and Daniel Litt just made available a preprint showing that there is a lower bound for the rank of geometric local systems with infinite monodromy on certain curves, and consequently the conjecture can not be true in this generality. 
}} 
\end{conj}

This density conjecture is in accordance with the one stated
in~\cite[Qu.~9.1]{EK19} and can be generalized to special loci, see Conjecture~\ref{conj:dens2}.
It is not difficult to show that Conjecture~\ref{dens:conj} holds for $G$ abelian. Indeed,
then  $G$ is a product of a unipotent group, a torus and a finite group, so one only has
to study the case that $G$ is a torus. However, the torus case follows from~\cite[Thm.~1.2]{EK19}.

A complex analytic  analog of the density conjecture involving the Riemann-Hilbert correspondence is
formulated in \cite[Conj.~10.4.1]{BW20}.

Another way to formulate Conjecture~\ref{dens:conj} is to
say that the image of
\[
\Ch_{G, \C}^\Box(\pi) \to \Ch_{\GL_{r,\C}}(\pi)
\]
contains a dense set of points corresponding to semi-simple representations $\pi\to
\GL_{r,\C}(\C)$  of geometric origin.
Here  \[\Ch_{\GL_{r,\C}}(\pi)=\Ch_{\GL_{r,\C}}^\Box(\pi)\sslash\GL_{r,\C}\] is the
character variety.

\medskip

We say that $\rho$ has \emph{quasi-unipotent
 monodromy at infinity} if for one (equivalently for all, see Proposition~\ref{prop:kas}) normal compactifications
$X\subset \overline X$
the
eigenvalues of
$\iota \circ\rho (T_D)$ are roots of unity.
Here $D\hookrightarrow
\overline{X}\setminus   X$ runs over the irreducible components which are of codimension one in $\overline X$ and $T_D$ is  the  canonical
conjugacy class $T_D\subset \pi$ corresponding to a  ``small loop around $D$'', see Section~\ref{sec:locmon}.

\begin{thm}[Monodromy theorem]
A $G$-representation $\rho\colon \pi \to G(\C)$ which is of geometric origin has
quasi-unipotent monodromy at infinity. 
\end{thm}

The monodromy theorem is due to Clemens and Landman,  see  \cite[Thm.~3.1]{Gri70}. Proofs
which are based on the study of local systems  were  given by Brieskorn
\cite[III,2]{Del70} and
Grothendieck \cite[Thm~1.2]{SGA7.1}. The proof of our main result,
Theorem~\ref{thm.main:intro} below, is motivated by Grothendieck's proof of the monodromy theorem.
In fact, in view of the monodromy theorem it can also be seen as a tiny bit of
evidence for the density conjecture.

\begin{thm}\label{thm.main:intro}[Theorem~\ref{thm.main1}]
The  set of $G$-representations  which have quasi-unipotent  monodromy at infinity is
Zariski dense in  $\Ch_{G,\C}^\Box(\pi) $.
\end{thm}

After we lectured on Theorem~\ref{thm.main:intro}, an alternative proof  for $G=\GL_{r,\C}$ based on the Riemann-Hilbert
correspondence and the Gelfond-Schneider theorem was given by B.~Bakker and Y.\
Brunebarbe.  Independently a similar density theorem involving the Riemann-Hilbert
correspondence was obtained by Budur, Lerer and Wang~\cite[Thm.~1.2]{BLW21}.

\medskip

Our  proof of  Theorem~\ref{thm.main:intro} is based  on the action of an  arithmetic
Galois group on certain completions of the character variety. This action is induced by
a comparison of the topological fundamental group with the \'etale fundamental group. On
the monodromy at infinity the Galois action is given in terms of the cyclotomic
character (see Lemma~\ref{lem.cyclchar}).

In Sections~\ref{sec:locmon} we recall some properties of fundamental groups. In
Section~\ref{sec:Grep} we introduce $G$-representations  of the fundamental group of a  variety with quasi-unipotent monodromy at infinity and
formulate our main theorem. The proof of the main theorem is contained in Section~\ref{sec:proofmain}.
In the final Section~\ref{sec:sploc} we explain how our proof can be applied more
generally to certain \emph{special} loci in $\Ch_{\GL_{r,\C}}^\Box(\pi)$ and how the 
density conjecture relates to the Fontaine-Mazur conjecture.

\medskip
{\em Acknowledgments:} We thank Pierre Deligne for a helpful discussion on the notion of
quasi-unipotent monodromy at infinity,  Ofer Gabber for mentioning the independence of $\iota$ for $\rho$ of geometric origin and Michel Brion for kindly supplying the reference \cite{Mil17} to us.

\section{The  monodromy at infinity}\label{sec:locmon}

\noindent \emph{The setup}.
The proof of our main result relies on arithmetic Galois groups, so we introduce a setting
in which we can later apply arithmetic arguments, even though the formulation of
Theorem~\ref{thm.main1} is purely complex.
Let $F\subset \mathbb C$ be a finitely generated field. Let $ X_0$ be a normal,
geometrically  irreducible variety over $F$ and let $  \overline X_0$ be a
normal compactification of $X_0$. Let $X\subset \overline X$ be the base change of theses
varieties to $\C$. We assume
that there exists a rational point $x_0\in X_0(F )$ and let $x$ be the associated
complex point. We consider the following objects:
\begin{itemize}
\item  $\bar F$ is the algebraic closure of $F$ in $\C$,  $\Gamma={\rm Gal}(\bar F/F)$ the
  Galois group of $F$.
  \item  $\pi=\pi_1(X(\C ),  x)$ is the topological fundamental group of $X(\C)$ based
    at $x$.
    \item
$\pi^\et=\pi_1^\et(X, x)\xrightarrow{\sim} \pi_1^{\rm
  \acute{e}t}(X_{0,\bar F}, x)$ is the geometric \'etale fundamental group of $X$  based at
$ x$, which by the Riemann existence theorem~\cite[Cor.~XII.5.2]{SGA1}  can be identified with the pro-finite completion of $\pi$.
\item $\pi_1^{\rm{\acute{e}t}} (X_0, x)$ is the arithmetic fundamental group of $X_0$ based at $ x$.
\item The conjugation action induced by the splitting
  of the homotopy exact  sequence~\cite[Thm.~IX.6.1] {SGA1}
\[
\xymatrix{
1\ar[r]  &  \pi^\et \ar[r]  &   \pi_1^\et(X_0, x) \ar[r] & \Gamma \ar[r]  \ar@{-->}@/_1.5pc/[l]& 1
}
\]
given by the point ${\rm Spec}(F) \xrightarrow{x_0}
  X \to {\rm Spec}(F)$ defines an action of $\Gamma$ on $\pi^{\rm \acute{e}t}$. 
  \item The family $D_i\hookrightarrow \overline X_0$ ($1\le i\le s$)  of irreducible components of
    $\overline X_0\setminus X_0$ which are of 
    codimension  one in $\overline X_0$. We assume that all $D_i$ are geometrically
    irreducible. By abuse of notation we denote the base change of $D_i$ to $\C$ by the
    same symbol.
\end{itemize}

\medskip

\noindent
\emph{Complex monodromy}.
To each $D_i$ one  associates a canonical conjugacy class $T_i\subset \pi$ as follows.
Consider the dense open subvariety
\[ \overline X^\circ=  \overline X\setminus  ( \overline X^{\rm sing}\cup (\overline X\setminus
  X)^{\rm sing} )\]
of $\overline X$.
Set $D^\circ_i=D_i\cap \overline X^\circ$ and $ X^\circ= \overline X^\circ \cap X$. Let us
assume $x\in  X^\circ(\C)$.
 Then  a ``small loop'' around $D^\circ_i(\C) \hookrightarrow X^\circ(\C)$ defines a canonical
 conjugation class  $T_i^\circ$ in $\pi_1(X^\circ(\C), x)$, see \cite[1.4]{Kas81}.  We define $T_i$ to be the image of $T_i^\circ$ via the surjective homomorphism 
$\pi_1(X^\circ(\C),x) \to \pi_1(X(\C),x)$.

\medskip

\noindent
\emph{\'Etale monodromy}.
We denote by $T^\et_i\subset \pi^\et$ the  conjugacy class induced by the image of $T_i\subset \pi$ in
$\pi^\et$.
It  can be described purely algebraically in
terms of  ramification theory, see \cite[XIV.1.1.10]{SGA7.2} for an exposition in the
one-dimensional case. This implies the following
well-known lemma, see also \cite[Claim~7.1]{EK20}.

\begin{lem}\label{lem.cyclchar}
For each $1\le i\le s$ the action of $\gamma\in \Gamma $ on $\pi^\et$ maps $T_i^\et$ to
$(T_i^\et)^{\chi(\gamma)}$.
Here $\chi\colon \Gamma\to \widehat{\mathbb Z}^\times$ is the cyclotomic character.
\end{lem}

\section{$G$-representations and quasi-unipotent monodromy} \label{sec:Grep}

\noindent \emph{Quasi-unipotent elements and $G$-representations}.
 Let $G/\C$ be a linear algebraic group. Recall that an element $g\in G(\C)$ is
called quasi-unipotent, if for one (or equivalently for any) embedding of algebraic groups $\iota\colon
G\hookrightarrow \GL_{r,\C}$ the eigenvalues of $\iota(g)$ are roots of unity.
Let $\pi$ be a finitely generated   group.
A \emph{ $G$-representation} of $\pi$  is a homomorphism 
$\rho\colon\pi\to G(\C)$.

 \medskip

 \noindent \emph{Character varieties}.
 Let $R $ be a noetherian ring and let $G$ be an affine group scheme of finite type $R$. There exists an affine scheme $\Ch_{G,R}^\Box(\pi)$  of
 finite type over $R$ such that for an $R$-algebra $R'$ there is a functorial bijection
\[
\Hom(\pi,G(R')) \cong \Ch_{G,R}^\Box(\pi)(R').
\]
The $R$-scheme $\Ch_{G,R}^\Box(\pi)$ is  called the \emph{framed character variety}. If
$\pi$ has a presentation $\langle w_1,\ldots , w_\sigma | r_1,\ldots ,r_\tau\rangle$ then 
\ga{1}{
\Ch_{G,R}^\Box(\pi) = \{ g \in    G^\sigma\,|\, r_1(g)=\cdots= r_\tau(g)=1 \}.
}

\medskip

\noindent
\emph{Quasi-unipotent monodromy at infinity}.
 Let the notation be as in Section~\ref{sec:locmon}, in particular $X$ is a normal complex
 variety, $x\in X(\C)$, $\pi=\pi_1(
 X(\C),x)$ and $X\subset \overline X$ is a normal compactification.
 Let $G/\C$ be a linear algebraic group.
We say that a $G$-representation $\rho\colon\pi\to G(\C)$ has \emph{quasi-unipotent monodromy at infinity} if
for all $1\le i\le s$ the image of the monodromy  $\rho(T_i)\subset G(\C )$ consists of
quasi-unipotent elements.
The following important theorem is shown in~\cite[Thm.3.1]{Kas81}. In fact Kashiwara's result 
is about constructible sheaves and one can easily
translate it into our setting of local systems.

\begin{prop}[Kashiwara] \label{prop:kas}
The property of $\rho\colon\pi\to G(\C)$ to have quasi-unipotent monodromy at infinity does not depend on
the choice of the normal compactification $\overline X$ of $X$.
\end{prop}

Our  main theorem says:

\begin{thm}\label{thm.main1}
The set of representations $\rho\in \Ch_{G,\C}^\Box(\pi)(\C)$ with quasi-unipotent monodromy at infinity
is Zariski dense in $ \Ch_{G,\C}^\Box(\pi)$.
\end{thm}

In Theorem~\ref{thm.main2} we  formulate a strengthening of Theorem~\ref{thm.main1}
involving an arithmetic Galois action.  Theorem~\ref{thm.main1} is shown in
Section~\ref{sec:proofmain}.

\begin{rmk}
One can also show by the same technique that the set of  representations $\rho\in
\Ch_{G,\C}^\Box(\pi)(\C)$ with finite determinant and with quasi-unipotent monodromy at
infinity is Zariski dense.
\end{rmk}

If $G$
is reductive we can form the categorical quotient of $\Ch_{G,K}(\pi) $  with respect to the conjugation action of $G$ to
obtain the \emph{character variety}
$$ \Ch_{G,\C}(\pi) =\Ch_{G,\C}^\Box(\pi)\sslash G .  $$
The points $\Ch_{G,K}(\pi)(\C)$
correspond to isomorphism classes of completely reducible representations $\rho\colon \pi \to G(\C)$  \cite[Sec.~11]{Sik12}. Theorem~\ref{thm.main1} then  immediately implies: 
\begin{cor}  \label{cor:main}
The set of isomorphism classes of completely reducible representations $\pi\to G(\C)$ with quasi-unipotent monodromy at infinity
is Zariski dense  in $\Ch_{G,\C}(\pi)$.
\end{cor}

\begin{ex}
For $X=\mathbb A^1\setminus (s \text{ points})$ the topological fundamental group $\pi=\pi_1(X(\C),x)$ is a free group with $s$
generators $w_1,\ldots, w_s$ (suitable loops around the $s$ points based at a common point $x\in
X(\C)$). The monodoromy at infinity for the canonical compactification $X\subset \P^1$
consists of the  conjugacy classes of $w_1,\ldots, w_s,(w_1\cdots w_s)^{-1}$ which
correspond to loops around the $s$ points $\A^1\setminus X$ and the point $\infty \in \P^1$. In
this case Theorem~\ref{thm.main1} says: The set of $g=(g_1,\ldots,g_s)\in G^s(\C)$ such that 
$g_1,\ldots, g_s,g_1\cdots g_s$ are quasi-unipotent is Zariski dense in $G^s$.

This example is related to \cite[Thm.\ B]{EK20} in the arithmetic situation. 
\end{ex}

\section{Proof of Theorem~\ref{thm.main1}}\label{sec:proofmain}

\noindent
\emph{$\Gamma$-action and  $\Ch_{G}^\Box(\pi)$}.
We use the notation of  Section~\ref{sec:locmon}, so $X$ is the base change  to $\C$ of a variety $X_0$ over a finitely generated field
$F\subset \C$.
Recall that $G\hookrightarrow \GL_{r,\C}$ is a linear algebraic group.

Let $Q\subset  \Ch_{G,\C}^\Box(\pi)$ be the Zariski closure of the set of quasi-unipotent
representations $\rho\colon \pi\to G(\C)$. We argue by contradiction and assume that $Q\ne
\Ch_{G,\C}^\Box(\pi)$. In particular,  $\Ch_{G,\C}^\Box(\pi)$ is non-empty.

Choose a  subring $R\subset \C$ which is of finite type over $\Z$, such that $G$ is
induced by a group scheme $\sG\hookrightarrow \GL_{r,R}$ over $R$ and such that  $Q$ is induced by a closed
   subscheme  $\sQ$ of $\Ch_{\sG,R}^\Box(\pi)$. Set $\sW=\Spec(R)$ and let $K\subset\C$
   be the field of fractions of $R$.

For a scheme $\sX$ of finite type over $R$, let us denote by $\vert \sX\vert $ the set of
closed points of $\sX$. For $x\in \sX$ we let $\sX^\wedge_x$ be the local scheme $\Spec(\sO_{\sX,x}^\wedge)$, where $\sO_{\sX,x}^\wedge$ is the completed local ring.

The $\Gamma$-action on $\pi^\et$ induces a continuous
$\Gamma$-action on the discrete set of closed points $\vert   \Ch_{\sG,R}^\Box(\pi) \vert$. Similarly, we get an induced $\Gamma_x$-action on $  \Ch_{\sG,R}^\Box(\pi)^\wedge_x$ for $x\in
\vert   \Ch_{\sG,R}^\Box(\pi) \vert$ and for $\Gamma_x\subset \Gamma$ the open stabilizer subgroup of $x$.

\medskip

\noindent
\emph{Characteristic polynomial of monodromy}.
For each local monodromy at infinity $T_i\subset \pi$, choose $g_i\in T_i$. 
  We have a morphism
   \ga{}{\psi\colon \Ch^\Box_{\sG,
       R}  \to \sN= \prod_{i=1}^s  (\A^{r-1}\times \G_m)  \notag }
 of affine schemes of finite type over $R$  
defined for each $i=1,\ldots, s$ by the coefficients
$(\sigma_1(\rho(g_i)),\ldots,\sigma_r(\rho(g_i)))\in \sN(R')$ of the characteristic polynomials
  \ga{}{     {\rm det}( T\cdot \mathbb I_r -\rho(g_i)) = T^r-\sigma_1(\rho(g_i)) T^{r-1} +
    \ldots + (-1)^r \sigma_r( \rho(g_i)) \notag}
  of a $G$-representation $\rho\colon \pi \to \sG(R')$, where $R'$ is an $R$-algebra.

  Furthermore, we have  the finite flat morphism
\ga{}{ \varphi\colon  \sM= (\G_m^r)^s\to \sN \notag} 
of affine schemes over $R$ given by 
\ga{}{ \G_m^r\to \A^{r-1}\times \G_m,  \ (\mu_1,\ldots,\mu_r)\mapsto (s_1(\mu_1,\ldots, \mu_r), \ldots, s_r(\mu_1, \ldots, \mu_r)), \notag}
where $s_i(\mu_1,\ldots, \mu_r)$ is the $i$-th elementary symmetric function in the
$\mu_j$.

The cyclotomic character $\chi$ induces an action of $\Gamma$ on $\vert \sM\vert$ and a
compatible action on $\vert\sN\vert$ such that $\vert \varphi\vert\colon \vert \sM\vert
\to \vert \sN\vert$ is $\Gamma$-equivariant. For each point $x\in \vert \sM\vert$ the stabilizer
$\Gamma_x\subset \Gamma$ acts on $\sM^\wedge_x$ and on $\sN^\wedge_x$ such that
$\varphi^\wedge_x$ is $\Gamma_x$-equivariant.

\medskip

\noindent
\emph{Certain closed points}.

  Let $\sT$ be the reduced closure  of the image of $\psi$. Let $\sS$ be
  $\varphi^{-1}(\sT)_{\rm red}$.
Note that the generic fibre $\sS_K$ of $\sS$ over $\sW$ is non-empty  as
$\Ch_{G,\C}^\Box(\pi)$   is non-empty,
 so  the smooth
locus $\sS^{\rm sm}$ of $\sS$ over $R$ is non-empty. By the generic flatness of $\psi$ we
can fix a closed point $z\in \Ch_{\sG,R}^\Box(\pi)\setminus \sQ$ such that
\begin{itemize}
\item $\psi$ is flat at $z$,
  \item $y=\psi(z)\in \varphi(\sS^{\rm sm})$.
\end{itemize}
We also fix a closed point $x\in\sS^{\rm sm}\cap \varphi^{-1}(y)$. Let $\Gamma'$ be the
intersection of stabilizers $\Gamma_x\cap\Gamma_z,$  which is thus open in $\Gamma$, and let $w\in \sW=\Spec(R)$ be the image of
the points $x,y,z$.

\begin{claim}
The closed subscheme $\sS^\wedge_x\hookrightarrow \sM^\wedge_x$ is $\Gamma'$-stable.
\end{claim}

\begin{proof}
  As $\psi$ is flat at the point $z$ the closed subscheme $\sT^\wedge_y\hookrightarrow
  \sN^\wedge_y $ is the schematic image of $\psi^\wedge_x\colon
  \Ch_{\sG,R}^\Box(\pi)^\wedge_z\to  \sN^\wedge_y$. As the latter morphism is
  $\Gamma'$-equivariant, it follows that $\sT^\wedge_y$ is stabilized by $\Gamma'$.
As $\sS^\wedge_x = \varphi^{-1}(\sT^\wedge_y)_{\rm red}$ we deduce that  $\sS^\wedge_x$  is stabilized by $\Gamma'$.
\end{proof}

\medskip

\noindent
\emph{De Jong's trick}.
For simplicity of notation we can assume that $\Gamma=\Gamma'$.
Choose a normal integral ring $A$ of finite type over $\Z$ with field of fractions $F$ such the characteristic of the residue field $k(x)$  of 
 $x$ is invertible in $A$. Then the action of $\Gamma$ on $\sM^\wedge_x$
via  the cyclotomic character $\chi$  factors
through   $\pi_1^\et(\Spec(A))$ and for an $\F_q$-point $a\colon\Spec(\F_q)\to \Spec(A)$
the associated Frobenius $\mathrm{Fr}=\mathrm{Fr}_a\in \pi_1^\et(\Spec(A))$, which is well-defined up to conjugation, acts
by multiplication by $q$ on the group scheme $\sM$ and on $\sM^\wedge_x$. 

\begin{claim}[De Jong's trick]\label{claim:dJtr}
  The morphism of local schemes
  \[
(\sS^\wedge_x)^{\rm Fr} \to \sW^\wedge_w
\]
is finite, flat and surjective.
\end{claim}

\begin{proof}
The following argument is copied   from~\cite[3.14]{deJ01}, see
also~\cite[Lem.~2.8]{Dri01}, ~\cite[Sec.~10]{EK19} and  ~\cite[Sec.~8]{EK20}. We can
assume without loss of generality that $k(x)=k(w)$.
By  smoothness of $\sS/\sW$ at $x$
\[
  \sS^\wedge_x \cong \Spec(\sO^\wedge_{\sW,w}\llbracket X_1,\ldots , X_j\rrbracket ).
\]
Then
\[
  (\sS^\wedge_x)^{\rm Fr} \cong \Spec ( \sO^\wedge_{\sW,w}\llbracket X_1,\ldots ,
  X_j\rrbracket/(1-\mathrm{Fr}(X_1), \ldots , 1-\mathrm{Fr}(X_j)  ) 
\]
has fibre dimension zero over $w$ as this fibre is a closed subscheme of the
$(q-1)$-torsion subscheme
of the torus $\sM_w$ over $w$. As in~\cite[3.14]{deJ01} basic commutative algebra shows that $(\sS^\wedge_x)^{\rm Fr}$ is  a local complete
intersection, finite and flat over $ \sW^\wedge_w$.
  \end{proof}

\medskip

\noindent
\emph{Conclusion}.
By Claim~\ref{claim:dJtr} there exists a point $\tilde x \in \sM$ which is $(q-1)$-torsion,
which maps to the generic point of $\sW$ and which specializes to $x$. In fact any point $\tilde x$  in the image of
the non-empty set
$(\sS^\wedge_x)^{\rm Fr}_K$ satisfies these properties. Then $\tilde y=\varphi(\tilde x)$
specializes to $y$. By going-down for flat morphisms there exists a point $\tilde z\in
\Ch_{\sG,R}^\Box(\pi)$ with $\varphi(\tilde z)=\tilde y$  which specializes to $z$. By
construction $\tilde z$ corresponds to a representation of $\pi$ which has quasi-unipotent
monodromy at infinity, so $\tilde z\in \sQ$ and therefore $z\in \sQ$. Contradiction!

\section{Special loci} \label{sec:sploc}

The aim of this section is to extend Conjecture~\ref{conj:dens2} and Theorem~\ref{thm.main1}
 to certain subloci of the character variety  $\Ch_{\GL_r,\C}(\pi)$. We also relate our
 density conjectures to other classical conjectures.  We use the notation of
Sections~\ref{sec:locmon}, \ref{sec:Grep} and~\ref{sec:proofmain}. Let $\wp \colon
\Ch_{\GL_r,\C}^\Box(\pi) \to \Ch_{\GL_r,\C}(\pi)$ be the canonical quotient map.
For a locally closed subscheme  $Z\hookrightarrow \Ch_{\GL_r,\C}(\pi)$, we denote by $\sZ\hookrightarrow
\Ch_{\GL_r,R}^\Box(\pi)$  a suitable locally closed subscheme such that $\wp(\sZ_{\mathbb C}
)\subset Z$ is dense. Here $R\subset \mathbb C$ is
a suitable subring of finite type over $\mathbb Z$ as above.

\begin{defn} \label{defn:sploc}
A  subscheme $Z\hookrightarrow \Ch_{\GL_r,\C}(\pi)$ as above is \emph{special} or
\emph{arithmetic} if there exist $R$ and $\sZ$ 
as above such that
for each closed point $z\in \sZ$, 
there is   an open subgroup of $\Gamma$ which stabilizes  the completion
$\sZ^\wedge_z$.
A point in $ s\in \Ch_{\GL_r,\C}(\pi)(\C)$ is \emph{ special}  or \emph{arithmetic} if the
subscheme $Z=\{ s \}$   is special. 
\end{defn}

 \begin{thm}\label{thm.main2}
 For a special subscheme $Z\hookrightarrow  \Ch_{\GL_r,\C}(\pi)$ the set of
 quasi-unipotent points in $Z$ is Zariski dense in $Z$. 
\end{thm}

\begin{proof}
The proof is analogous to the one of Theorem~\ref{thm.main1}. One just replaces
$\Ch_{G,\C}^\Box(\pi)$ by $\wp^{-1}(Z)$ and $Q$ by $\wp^{-1}(Z)\cap Q$.
\end{proof}

Here is another natural density conjecture in this context:

\begin{conj}[Density]\label{conj:dens2}
Let  $Z\hookrightarrow
\Ch_{ \GL_r,\C}(\pi)$ be a special subscheme. Then the set of complex points of $Z$
corresponding of  to representations of geometric origin $\rho\colon \pi\to \GL_r (\mathbb
C)$ is dense.

In particular the special points are then dense on $Z$. {\footnote{See the footnote to Conjecture~\ref{dens:conj}}.}
\end{conj}

Note that a point of $  \Ch_{\GL_r,\C}(\pi)(\C)$ which corresponds to a representation of
geometric origin  is special by the comparison isomorphism between Betti cohomology and
$\ell$-adic cohomology.

\smallskip

The following observations are easy to check.

\begin{rmk}\label{rmk.special}\mbox{}
  \begin{itemize}
  \item[(1)]
    Conjecture~\ref{conj:dens2} $ \;\Rightarrow\; $  Conjecture~\ref{dens:conj}.   Indeed
    for $G \stackrel{\iota}{\hookrightarrow} \GL_{r,\C }$ given, the image of  $ \Ch_{G,\C}^\Box(\pi) \to 
\Ch_{ \GL_r,\C}(\pi)$ is construcible and we take $Z$ to be a suitable dense subscheme in
this image.

    \item[(2)]
 Conjecture~\ref{conj:dens2} comprises the density conjecture  formulated
 in~\cite[Qu.~9.1]{EK19}   1), 2), except 3). For 1) and 2) this is by definition.
 \item[(3)]  Conjecture~\ref{conj:dens2} for $\dim(Z)=0$ implies Simpson's
``rigid $\Rightarrow$ motivic'' conjecture~\cite[Conj.~4]{Sim90} as 
rigid representations are arithmetic~\cite[Thm.~4]{Sim92}.
\item[(4)] Recent work of Petrov~\cite{Pet20} implies that the relative Fontaine--Mazur
  conjecture of Liu--Zhu~\cite{LZ17} implies  Conjecture~\ref{conj:dens2} for $\dim(Z)=0$.
\end{itemize}
\end{rmk}

In fact (4) is a direct consequence of~\cite[Lem.~6.2]{Pet20} over a number field $F$. For
a general finitely generated field $F\subset \mathbb C$ one has to use a spreading 
argument similar to \cite[Prop.~6.1]{Pet20} in order to reduce to the number field case. 

So (3) and (4) of Remark~\ref{rmk.special} together say that

\[
  \boxed{
 \text{ relative F-M conj.} \;\Rightarrow\;  \text{ Simpson's
  ``rigid $\Rightarrow$ motivic'' conj. }}
  \]

\end{document}